\documentclass[a4paper,10pt]{article}

\usepackage[latin1]{inputenc}
\usepackage{latexsym,amssymb,amsmath,epsfig}
\usepackage{a4wide}
\usepackage{graphicx}
\usepackage{theorem}
\theorembodyfont{\upshape}

\def\RR{\mbox{\it I\hskip -0.177em R}}
\def\PP{\mbox{\it I\hskip -0.177em P}}
\def\NN{\mbox{\it I\hskip -0.177em N}}

\def\EE{\mbox{\it I\hskip -0.177em E}}

\newcommand{\sRR}{{\hbox{$\scriptstyle{I}$\kern-.25em\hbox
{$\scriptstyle{R}$}}}}
\newcommand{\sPP}{{\hbox{$\scriptstyle{I}$\kern-.25em\hbox
{$\scriptstyle{P}$}}}}
\newcommand{\sNN}{{\hbox{$\scriptstyle{I}$\kern-.25em\hbox
{$\scriptstyle{N}$}}}}
\newcommand{\sZZ}{{\hbox{$\scriptstyle{Z}$\kern-.3em\hbox
{$\scriptstyle{Z}$}}}}
\newcommand{\sEE}{{\hbox{$\scriptstyle{I}$\kern-.25em\hbox
{$\scriptstyle{E}$}}}}
\newcommand{\sII}{{\hbox{$\scriptstyle{I}$\kern-.25em\hbox
{$\scriptstyle{I}$}}}}

\newtheorem{defi}{Definition}[section]
\newtheorem{lemm}[defi]{Lemma}
\newtheorem{prop}[defi]{Proposition}

\newtheorem{theo}[defi]{Theorem}
\newenvironment{proof}{\vskip 2mm\noindent {\it Proof}:}
                    {\hfill $\square$ \vskip 2mm \noindent}
\newenvironment{proofbis}[1]{\vskip 2mm\noindent {\it Proof (of #1)}:}
                    {\hfill $\square$ \vskip 2mm \noindent} 

\title{Poisson approximations for the Ising model}
\author{D. Coupier} 
\date{\today}

\begin{document}
\maketitle

\begin{center}
Ecole Polytechnique Universitaire de Lille
\end{center}
\vskip 0.5cm
{\bf E-Mail address :} {\tt david.coupier@polytech-lille.fr}\\  

\vskip 0.3cm 
\noindent
{\bf Mail address :} Ecole Polytechnique Universitaire de Lille,\\
\hspace*{2.4cm} Cit\'e scientifique - Avenue Paul Langevin\\
\hspace*{2.4cm} 59655 Villeneuve d'Ascq Cedex - France.

\noindent
{\bf Telephone :} 33 3 28 76 74 36 \\
{\bf Fax :} 33 3 28 76 73 21
\vskip 2cm

\begin{abstract}
A $d$-dimensional Ising model on a lattice torus is considered. As the size $n$ of the lattice tends to infinity, a Poisson approximation is given for the distribution of the number of copies in the lattice of any given local configuration, provided the magnetic field $a=a(n)$ tends to $-\infty$ and the pair potential $b$ remains fixed. Using the Stein-Chen method, a bound is given for the total variation error in the ferromagnetic case.
\end{abstract}
\vskip 1cm

\noindent {\bf Key words :} Poisson approximation, Ising model, ferromagnetic interaction, Stein-Chen method.
\vskip 0.5cm
\noindent
{\bf AMS Subject Classification :} 60F05, 82B20.
\newpage

\section{Introduction}
\label{section:intro}

The following situation, called ``the law of small numbers'', is very classical in probability theory. Suppose $\{I_{\lambda}\}_{\lambda\in\Lambda}$ is a finite family of indicator random variables, with the properties that the probabilities $\PP(I_{\lambda}=1)$ are small and that there is not too much dependence between the $I_{\lambda}$'s. Then, it is reasonable to expect the distribution of $\sum_{\lambda\in\Lambda} I_{\lambda}$ to be approximately Poisson. In the theory of random graphs, inaugurated by Erd{\"o}s and R\'enyi \cite{ErdosRenyi}, such results are frequent (see \cite{Bollobas} or \cite{Spencer} for a general reference). The $I_{\lambda}$ can, for instance, indicate the places in the random graph where a given subgraph appear. Some analogous results hold for random colorings of a lattice graph in dimension $2$, corresponding to the context of random images \cite{CoupierDesolneuxYcart2}. In both cases, the models are built on a large number of independent random variables: the edges of a random graph or the pixels of a random image. In this article, we shall study Poisson approximations for sums of indicators defined from a large number of dependent random variables, namely the spins of an Ising model.\\
Let us consider a lattice graph in dimension $d\geq 1$, with periodic boundary conditions (lattice torus). The vertex set is $V_{n}=\{ 0,\ldots,n-1\}^{d}$. The integer $n$ will be called the \textit{size} of the lattice. The edge set, denoted by $E_{n}$, will be specified by defining the set of neighbors $\mathcal{V}(x)$ of a given vertex $x$:
\begin{equation}
\label{def:voisinage}
\mathcal{V}(x) = \{ y \neq x \in V_{n}\,,\; \|y-x\|_{p} \leq \rho\} ~,
\end{equation}
where the substraction is taken componentwise modulo $n$, $\| \cdot \|_{p}$ stands for the $L_{p}$ norm in $\mathbb{R}^{d}$ ($1\leq p\leq\infty$), and $\rho$ is a fixed parameter. For instance, the square lattice is obtained for $p=\rho=1$. Replacing the $L_{1}$ norm by the $L_{\infty}$ norm adds the diagonals. From now on, all operations on vertices will be understood modulo $n$. In particular, each vertex of the lattice has the same number of neighbors.\\
A \textit{configuration} is a mapping from the vertex set $V_{n}$ to the state space $W = \{ -1, +1 \}$. Their set is denoted by $\mathcal{X}_{n} = W^{V_{n}}$ and called the \textit{configuration set}. Here, we shall deal with one of the simplest and most widely studied parametric families of random field distributions: the Ising model (see e.g. \cite{Georgii,Malyshev}).

\begin{defi}
\label{def:ising}
\textit{
Let $G_{n}=(V_{n},E_{n})$ be an undirected graph structure with finite vertex set $V_{n}$ and edge set $E_{n}$. Let $a$ and $b$ be two reals. The \textnormal{Ising model} with parameters $a$ and $b$ is the probability measure $\mu_{a,b}$ on $\mathcal{X}_{n}=\{ -1,+1\}^{V_{n}}$ defined by: $\forall\sigma\in\mathcal{X}_{n}$,
\begin{equation}
\label{def:mes_gibbs}
\mu_{a,b} (\sigma) = \frac{1}{Z_{a,b}} \exp \left( a \sum_{x \in V_{n}} \sigma(x) + b \sum_{ \{x,y\} \in E_{n}} \sigma(x) \sigma(y) \right)\;,
\end{equation}
where the normalizing constant $Z_{a,b}$ is such that $\sum_{\sigma\in \mathcal{X}_{n}} \mu_{a,b}(\sigma)=1$.
}
\end{defi}

Following the definition of \cite{Malyshev} p.~2, the measure $\mu_{a,b}$ defined above is a Gibbs measure associated to potentials $a$ and $b$. Expectations relative to $\mu_{a,b}$ will be denoted by $\EE_{a,b}$.\\
In the classical presentation of statistical physics, the elements of $\mathcal{X}_{n}$ are spin configurations; each vertex of $V_{n}$ is an atom whose spin is either positive or negative. Here, we shall simply talk about positive or negative vertices instead of positive or negative spins. The parameters $a$ and $b$ are respectively the \textit{magnetic field} and the \textit{pair potential}. The model remaining unchanged by swapping positive and negative vertices and replacing $a$ by $-a$, we chose to study only negative values of the magnetic field $a$.\\
Various ``laws of small numbers'' have been already proved for the Ising model. Fern\'andez et al. \cite{FFG1,FFG2} have established the asymptotic Poisson distribution of contours in the nearest-neighbor Ising model at low temperature (i.e. $b$ large enough) and zero magnetic field. The Stein-Chen method is a useful way to get Poisson approximations; see \cite{AGG}, \cite{Barbour} for a very complete reference or \cite{Chen} for the original paper of Chen. Barbour and Greenwood \cite{BarbourGreenwood} have applied it to a class of Markov random fields; the bounds that they obtained for the Ising model are not quite explicit. In the same context as \cite{FFG2}, Ferrari and Picco \cite{FerrariPicco} have found bounds on the total variation distance between the law of large contours and a Poisson process. Ganesh et al. \cite{GHOSU} have studied the Ising model for positive values of $b$. Provided the magnetic field $a$ tends to infinity, they proved that the distribution of the number of negative vertices is approximately Poisson.\\
Our goal is to generalize the convergence in distribution given by \cite{GHOSU} to any value of the pair potential $b$ and to objects more elaborated than a single vertex.\\
We are interested in the occurrences in the graph $G_{n}$ of a fixed \textit{local configuration} $\eta$ (see Section \ref{section:proba_cond} for a precise definition and Figure \ref{fig:local_configuration} for an example). Such a configuration is called ``local'' in the sense that the vertex set on which it is defined is fixed and does not depend on $n$. Its number of occurrences in $G_{n}$ is denoted by $X_{n}(\eta)$.\\
As the size $n$ of the lattice tends to infinity, the potential $a=a(n)$ will depend on $n$ whereas the potential $b$ will remain fixed. The case where $a(n)$ tends to $-\infty$ corresponds to rare positive vertices among a majority of negative ones. As a consequence, the local configuration $\eta$ may occur or not in the graph, depending on its number of positive vertices $k(\eta)$. See Proposition $4.2$ of \cite{CoupierDoukhanYcart} for a precise description of this phenomenon. In particular, in order to get a nontrivial limiting result for the probability $\mu_{a,b}(X_{n}(\eta)>0)$, it is needed to take $e^{2a(n)}$ of order $n^{-d/k(\eta)}$. Therefore, throughout this paper, the magnetic field $a(n)$ will satisfy the identity
\begin{equation}
\label{ordre_a(n)}
e^{2 a(n)} = c n^{-d/k(\eta)} ~,
\end{equation}
where $c$ is a positive constant. Our first result describes the asymptotic behavior of the number $X_{n}(\eta)$ of occurrences of $\eta$ in the lattice: it will be poissonian and will depend on $k(\eta)$ through (\ref{ordre_a(n)}), but also on the geometry of $\eta$ through its \textit{perimeter}, denoted by $\gamma(\eta)$.

\begin{theo}
\label{theo:PA}
\textit{
Assume that the magnetic field $a(n)$ satisfies \textnormal{(\ref{ordre_a(n)})} and that the pair potential $b$ is an arbitrary real number. As $n$ tends to infinity, the distribution of $X_{n}(\eta)$ converges weakly to the Poisson distribution with parameter $c^{k(\eta)} e^{-2 b \gamma(\eta)}$.
}
\end{theo}

The proof is based on the moment method (see \cite{Bollobas} p.~25 or Lemma \ref{lemm:moment}), and requires estimates based on the \textit{local energy} of $\eta$ (Definition \ref{def:energy}). The result of Ganesh et al. \cite{GHOSU} is obtained as a particular case when the pair potential $b$ is positive and $\eta$ is a single positive vertex: $k(\eta)=1$ and $\gamma(\eta)=4$.\\
The Stein-Chen method makes it possible to obtain good estimates on the accuracy of Poisson approximations in terms of \textit{total variation distance}. When the Gibbs measure $\mu_{a,b}$ defined in (\ref{def:mes_gibbs}) satisfies the FKG inequality \cite{FKG} (i.e. for positive values of the pair potential $b$), this method is applied to a sum of \textit{increasing} random indicators (\ref{sum_increasing}) and produces Lemma \ref{lemm:chenstein}. Then, bounds on the first two moments of the random variable $X_{n}(\eta)$ (Lemmas \ref{lemm:T1} and \ref{lemm:M2}) allow to precise the Poisson approximation given by Theorem \ref{theo:PA}. This leads to our second result, where $\mathcal{L}(X)$ and $\mathcal{P}(\lambda)$ respectively denote the distribution of $X$ and the Poisson distribution with parameter $\lambda$.

\begin{theo}
\label{theo:TV}
\textit{
Assume that the magnetic field $a(n)$ satisfies \textnormal{(\ref{ordre_a(n)})} and that the pair potential $b$ is positive. Then, the total variation distance between $\mathcal{L}(X_{n}(\eta))$ and the Poisson distribution with parameter $c^{k(\eta)} e^{-2b \gamma(\eta)}$ satisfies:
$$
d_{TV}(\mathcal{L}(X_{n}(\eta)), \mathcal{P}(c^{k(\eta)} e^{-2b \gamma(\eta)})) = \mathcal{O}(n^{-d/k(\eta)}) ~.
$$
}
\end{theo}

The paper is organized as follows. The notion of local configuration $\eta$ is defined in Section \ref{section:proba_cond}. Its number of positive vertices $k(\eta)$ and its perimeter $\gamma(\eta)$ are also introduced. Lemma \ref{lemm:propre} reduces proofs of Theorems \ref{theo:PA} and \ref{theo:TV} to \textit{clean} local configurations. In this case, integers $k(\eta)$ and $\gamma(\eta)$ naturally occur in the expression of the local energy of $\eta$. Describing this quantity will be essential in our study. This allows us to control the conditional probability of $\eta$ to occur in the graph (Lemma \ref{lemm:proba_cond}). It immediatly follows that the expected number of occurrences of $\eta$ in $G_{n}$ tends to $c^{k(\eta)} e^{-2b \gamma(\eta)}$, as $n$ tends to infinity. Finally, Sections \ref{section:PA} and \ref{section:ferromagnetic} are respectively devoted to the proofs of Theorems \ref{theo:PA} and \ref{theo:TV}.

\section{Conditional probability of a local configuration}
\label{section:proba_cond}


Let us start with some notations and definitions. Given $\sigma\in \mathcal{X}_{n}=W^{V_{n}}$ and $V\subset V_{n}$, we denote by $\sigma_{V}$ the natural projection of $\sigma$ over $W^{V}$. If $U$ and $V$ are two disjoint subsets of $V_{n}$ then $\sigma_{U}\sigma'_{V}$ is the configuration on $U\cup V$ which is equal to $\sigma$ on $U$ and $\sigma'$ on $V$. Let us denote by $\delta V$ the neighborhood of $V$ (corresponding to (\ref{def:voisinage})):
$$
\delta V = \{ y \in V_{n} \setminus V , \; \exists x \in V , \; \{ x,y \} \in E_{n} \} ~,
$$
and by $\overline{V}$ the union of the two disjoint sets $V$ and $\delta V$. Moreover, $|V|$ denotes the cardinality of $V$ and $\mathcal{F}(V)$ the $\sigma$-algebra generated by the configurations of $W^{V}$.\\
As usual, the graph distance $dist$ is defined as the minimal length of a path between two vertices. We shall denote by $B(x,r)$ the ball of center $x$ and radius $r$:
$$
B(x,r) = \{\, y\in V_{n}\,;\; dist(x,y)\leq r\,\} ~.
$$
In the case of balls, $\overline{B(x,r)}=B(x,r+1)$. In order to avoid unpleasant situations, like self-overlapping balls, we will always assume that $n>2\rho r$. If $n$ and $n'$ are both larger than $2\rho r$, the balls $B(x,r)$ in $G_n$ and $G_{n'}$ are isomorphic. Two properties of the balls $B(x,r)$ will be crucial in what follows. The first one is that two balls with the same radius are translates of each other:
$$
B(x+y,r)=y+B(x,r)\;.
$$
The second one is that for $n>2\rho r$, the cardinality of $B(x,r)$ depends only on $r$ and neither on $x$ nor on $n$: it will be denoted by $\beta(r)$. The same is true for the number of edges $\{y,z\}\in E_{n}$ with $y,z\in B(x,r)$, which will be denoted by $\alpha(r)$.\\

\begin{figure}[!ht]
\begin{center}
\includegraphics[width=4.5cm,height=4.5cm]{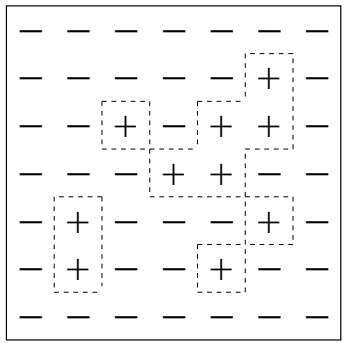}
\end{center}
\caption{\label{fig:local_configuration}A clean local configuration $\eta$ with $k(\eta)=|V_{+}(\eta)|=10$ positive vertices, in dimension $d=2$ and on a ball of radius $r=3$ (with $\rho=1$ and relative to $\| \cdot \|_{\infty}$).}
\end{figure}

Let $r$ be a positive real, and consider a fixed ball with radius $r$, say $B(0,r)$. We denote by $\mathcal{D}_{r}=W^{B(0,r)}$ the set of configurations on that ball. Elements of $\mathcal{D}_{r}$ will be called \textit{local configurations of radius $r$}. A local configuration $\eta\in\mathcal{D}_{r}$ is determined by its subset $V_{+}(\eta)\subset B(0,r)$ of positive vertices:
$$
V_{+}(\eta) = \{ x \in B(0,r), \; \eta(x) = +1 \} ~.
$$
The cardinality of this set will be denoted by $k(\eta)$ and its complementary set in $B(0,r)$, i.e. the set of negative vertices of $\eta$, by $V_{-}(\eta)$. Of course, there exists only a finite number of local configurations of radius $r$ (precisely $2^{\beta(r)}$). In what follows, $\eta$, $\eta'$ will denote local configurations of radius $r$ and $\zeta$, $\zeta'$ those of radius larger than $r$.\\
A local configuration $\eta\in\mathcal{D}_{r}$ is said \textit{clean} if its subset of plus vertices $V_{+}(\eta)$ is included in the ball $B(0,r-1)$. In other words, vertices of a clean local configuration which are at distance $r$ from the center $0$ are negative. Figure \ref{fig:local_configuration} shows an example of such a local configuration.

Let $\eta\in\mathcal{D}_{r}$. For each vertex $x\in V_{n}$, denote by $\eta_{x}$ the translation of $\eta$ onto the ball $B(x,r)$ (up to periodic boundary conditions): 
$$
\forall y \in V_{n}, \; dist(0,y)\leq r \Longrightarrow \eta_{x}(x+y) = \eta(y) ~.
$$
Let us denote by $I_{x}^{\eta}$ the indicator function defined on $\mathcal{X}_{n}$ as follows: $I_{x}^{\eta}(\sigma)$ is $1$ if the restriction of the configuration $\sigma\in \mathcal{X}_{n}$ to the ball $B(x,r)$ is $\eta_{x}$ and $0$ otherwise. Finally, let us define the random variable $X_{n}(\eta)$ which counts the number of copies of the local configuration $\eta$ in $G_{n}$:
$$
X_{n}(\eta) = \sum_{x \in V_{n}} \; I_{x}^{\eta} ~.
$$
Due to periodicity, this sum consists of $n^{d}$ indicator functions $I_{x}^{\eta}$, which have the same distribution.\\


In order to control the random variable $X_{n}(\eta)$, we describe its ``local behavior'' by introducing the \textit{local energy} of $\eta$. Let us start with the following definition.

\begin{defi}
\label{def:energy}
\textit{
Let $x\in V_{n}$ and $\sigma\in W^{B(x,r+1)}$. The \textnormal{local energy} $H^{B(x,r)}(\sigma)$ of the configu\-ration $\sigma$ on the ball $B(x,r)$ is defined by:
$$
H^{B(x,r)}(\sigma) = a(n) \sum_{y \in B(x,r)} \sigma(y) + b \sum_{{\scriptstyle \{ y,z \} \in E_{n}} \atop \scriptstyle (y \in B(x,r)) \vee (z \in B(x,r))} \sigma(y) \sigma(z) ~,
$$
where $(y\in B(x,r))\vee (z\in B(x,r))$ means at least one of the two vertices $y$ and $z$ belongs to $B(x,r)$ (the other might belong to its neighborhood $\delta B(x,r)$).
}
\end{defi}

Let us fix a vertex $x$ and denote merely by $B$ the ball $B(x,r)$. For any local configuration $\eta\in\mathcal{D}_{r}$ and for any $\sigma\in W^{\delta B}$, the local energy $H^{B}(\eta_{x} \sigma)$ on $B$ of the configuration which is equal to $\eta_{x}$ on $B$ and $\sigma$ on $\delta\!B$ can be expressed as:
\begin{equation}
\label{energie_explicite}
a(n) (2k(\eta) - \beta(r)) + b \left( \sum_{{\scriptstyle \{ y,z \} \in E_{n}} \atop \scriptstyle y,z \in B} \eta_{x}(y) \eta_{x}(z) + \sum_{{\scriptstyle \{ y,z \} \in E_{n}} \atop \scriptstyle y \in B, z \in \delta B} \eta_{x}(y) \sigma(z) \right) ~.
\end{equation}
Actually, this notion of local energy allows us to explicitly write the conditional probability $\mu_{a,b}(I_{x}^{\eta}=1 | \sigma)$, $\sigma \in W^{\delta\!B}$:
\begin{equation}
\label{proba_cond}
\mu_{a,b}(I_{x}^{\eta}=1 | \sigma) = \frac{e^{H^{B}(\eta_{x} \sigma)}}{\sum_{\eta' \in \mathcal{D}_{r}} e^{H^{B}(\eta_{x}' \sigma)}} ~.
\end{equation}
As we shall see in Lemma \ref{lemm:proba_cond}, bounding the above conditional probability is central in our study.\\
An easy way to connect the number of copies of the local configuration $\eta$ to its local energy consists in writing, for any given vertex $x$:
\begin{eqnarray}
\label{egal_expectation}
\EE_{a,b} \lbrack X_{n}(\eta) \rbrack & = & \EE_{a,b} \lbrack n^{d} I_{x}^{\eta} \rbrack \nonumber \\
& = & \EE_{a,b} \lbrack n^{d} \mu_{a,b}(I_{x}^{\eta}=1 | \mathcal{F}(\delta B)) \rbrack ~.
\end{eqnarray}
Here, $\mu_{a,b}(I_{x}^{\eta}=1 | \mathcal{F}(\delta\!B))$ represents a $\mathcal{F}(\delta\!B)$-measurable random variable and, for $\sigma \in W^{\delta\!B}$, $\mu_{a,b}(I_{x}^{\eta}=1 | \mathcal{F}(\delta\!B))(\sigma) = \mu_{a,b}(I_{x}^{\eta}=1 | \sigma)$ a conditional probability. Note that the set $\delta\!B$ has bounded cardinality (not depending on $n$). Then, from a convergence result on the random variable $n^{d} \mu_{a,b}(I_{x}^{\eta}=1 | \mathcal{F}(\delta\!B))$ it will be easy to obtain a similar result for its expectation, i.e. for $\EE_{a,b} \lbrack X_{n}(\eta) \rbrack$.\\


We will now give the reason for the hypothesis (\ref{ordre_a(n)}), that links the magnetic field $a(n)$ to the number of positive vertices of the local configuration $\eta$. The event $X_{n}(\eta)>0$ corresponds to the appearance of $\eta$ in the graph $G_{n}$. In \cite{CoupierDoukhanYcart} Proposition $4.2$, it has been proved that:
$$
\mbox{if } \lim_{n\to\infty} e^{2 a(n) k(\eta)} n^{d} = 0 \mbox{ then }
$$
\begin{equation}
\label{threshold1}
\lim_{n\to\infty} \mu_{a,b}(X_{n}(\eta)>0) = 0 \mbox{ and } \lim_{n\to\infty} \EE_{a,b}\lbrack X_{n}(\eta) \rbrack = 0 ~;
\end{equation}
$$
\mbox{if } \lim_{n\to\infty} e^{2 a(n) k(\eta)} n^{d} = +\infty \mbox{ then }
$$
\begin{equation}
\label{threshold2}
\lim_{n\to\infty} \mu_{a,b}(X_{n}(\eta)>0) = 1 \mbox{ and } \lim_{n\to\infty} \EE_{a,b}\lbrack X_{n}(\eta) \rbrack = +\infty ~.
\end{equation}
\noindent
In particular, the element of $\mathcal{D}_{r}$ having only negative vertices, called the \textit{null local configuration} and denoted by $\eta^{0}$, has a probability which always tends to $1$.\\
From now on, assume that $\eta$ has at least one positive vertex; $k(\eta)\geq 1$. Using the vocabulary of the random graph theory, statements (\ref{threshold1}) and (\ref{threshold2}) essentially mean that the quantity $n^{-d/k(\eta)}$ is the \textit{threshold function} (for $e^{2a(n)}$) of the property $X_{n}(\eta)>0$. It does not depend on the radius $r$ of the ball on which the local configuration $\eta$ is defined: $r$ is just a phantom parameter which serves only to ensure that $\eta$ is a local configuration. Actually, the function $n^{-d/k(\eta)}$ only depends on the number of positive vertices of $\eta$. Roughly speaking, if $e^{2a(n)}$ is small compared to $n^{-d/k(\eta)}$, then asymptotically, there is no copy of $\eta$ in $G_{n}$. If $e^{2a(n)}$ is large compared to $n^{-d/k(\eta)}$, then at least one copy of $\eta$ can be found in the graph, with probability tending to $1$.\\
Consequently, in order to get a nontrivial limiting result for the probability $\mu_{a,b}(X_{n}(\eta)>0)$, it is needed to take $e^{2a(n)}$ of order $n^{-d/k(\eta)}$. Hence, for the rest of this article, (\ref{ordre_a(n)}) is satisfied, i.e.
$$
e^{2a(n)} = c n^{-d/k(\eta)} ~,
$$
for some positive constant $c$. Under this hypothesis, statement (\ref{threshold1}) says that asymptotically there will be no local configurations with (strictly) more than $k(\eta)$ positive vertices in the lattice. The following lemma quantifies this result.

\begin{lemm}
\label{lemm:fonda}
\textit{
Let $\eta\in\mathcal{D}_{r}$ and suppose that $e^{2a(n)}=cn^{-d/k(\eta)}$, for some constant $c>0$. Let $R\geq r$ an integer and $\zeta\in\mathcal{D}_{R}$. Then, there exists a constant $M_{1}>0$ such that for all $n$, for all vertex $x\in V_{n}$ and for all configuration $\sigma\in\mathcal{X}_{n}$,
\begin{equation}
\label{major_zeta}
n^{d} \mu_{a,b}(I_{x}^{\zeta}=1 | \sigma_{\delta B(x,R)}) \leq M_{1} e^{2 a(n)(k(\zeta) - k(\eta))} ~,
\end{equation}
and for all $n$,
\begin{equation}
\label{major_esperance_zeta}
\EE_{a,b} \lbrack X_{n}(\zeta) \rbrack \leq M_{1} e^{2 a(n)(k(\zeta) - k(\eta))} ~.
\end{equation}
}
\end{lemm}

\begin{proof}
Let $x\in V_{n}$ and denote merely by $B_{R}$ the ball $B(x,R)$. For any configuration $\sigma\in W^{\delta B_{R}}$, the local energy $H^{B_{R}}(\zeta_{x}\sigma)$ is given by (\ref{energie_explicite}). Since the set $\overline{B}_{R}$ has a bounded cardinality, there exists a constant $K>0$, only depending on the pair potential $b$ and the radius $R$ (and not on $n$, nor $x$, nor $\sigma$), such that
$$
H^{B_{R}}(\zeta_{x} \sigma) - H^{B_{R}}(\zeta^{0}_{x} \sigma) \leq 2 a(n) k(\zeta) + K ~,
$$
where $\zeta^{0}$ is the null local configuration of radius $R$ ($k(\zeta^{0})=0$). Hence, by relation (\ref{proba_cond}) and for all $\sigma\in W^{\delta B_{R}}$:
\begin{eqnarray*}
\mu_{a,b}(I_{x}^{\zeta} = 1 | \sigma ) & \leq & e^{H^{B_{R}}(\zeta_{x} \sigma) - H^{B_{R}}(\zeta^{0}_{x} \sigma)} \\
& \leq & e^{2 a(n) k(\zeta) + K} ~.
\end{eqnarray*}
Finally, hypothesis (\ref{ordre_a(n)}) provides the first inequality of Lemma \ref{lemm:fonda}:
\begin{eqnarray*}
n^{d} \mu_{a,b}(I_{x}^{\zeta}=1 | \sigma ) & \leq & n^{d} e^{2 a(n) k(\eta) + K} e^{2 a(n)(k(\zeta) - k(\eta))} \\
& \leq & c^{k(\eta)} e^{K} e^{2 a(n)(k(\zeta) - k(\eta))} ~,
\end{eqnarray*}
with $M_{1}=c^{k(\eta)}e^{K}>0$. The quantity $n^{d}\mu_{a,b}(I_{x}^{\zeta}=1|\sigma)$ is bounded uniformly on the configuration $\sigma\in W^{\delta B_{R}}$. So, its expectation satisfies the same inequality and (\ref{major_esperance_zeta}) follows.
\end{proof}

A primary consequence of Lemma \ref{lemm:fonda} consists in reducing our study to clean local configurations. To any given $\eta\in\mathcal{D}_{r}$, a local configuration $\mathring{\eta}\in\mathcal{D}_{r+1}$ is associated by the following process:
$$
\mathring{\eta}(x) = \left\lbrace \begin{array}{l}
\eta(x) \; \mbox{ if } \; x \in B(0,r) ~,\\
-1 \; \mbox{ if } \; dist(x,0) = r+1 ~.
\end{array} \right.
$$
The local configuration $\mathring{\eta}$ is clean and satisfies $k(\mathring{\eta})=k(\eta)$. Note that the inequality $X_{n}(\mathring{\eta})\leq X_{n}(\eta)$ holds for all size $n$. Actually, these two random variables are asymptotically equal. Indeed, assume that $\eta$ occurs on the ball $B(x,r)$. Then, hypothesis (\ref{ordre_a(n)}) forces vertices at distance $r+1$ from $x$ to be negative with probability tending to $1$:
$$
\lim_{n\to+\infty} \mu_{a,b} \left( I_{x}^{\mathring{\eta}} = 1 | I_{x}^{\eta} = 1 \right) = 1 ~.
$$
Lemma \ref{lemm:propre} expresses this result in terms of \textit{total variation distance}. Recall that if $\mu$ and $\nu$ are two probability distributions, the total variation distance between $\mu$ and $\nu$ is
$$
d_{TV}(\mu, \nu) = \sup_{A} | \mu(A) - \nu(A) | ~,
$$
where the supremum is taken over all measurable sets. Besides, the probability distribution of a random variable $X$ will be denoted by $\mathcal{L}(X)$.

\begin{lemm}
\label{lemm:propre}
\textit{
Let $\eta\in\mathcal{D}_{r}$ and suppose that $e^{2a(n)}=cn^{-d/k(\eta)}$, for some constant $c>0$. Then, the total variation distance between distributions of $X_{n}(\mathring{\eta})$ and $X_{n}(\eta)$ satisfies:
\begin{equation}
\label{dvt_propre}
d_{TV}\left( \mathcal{L}(X_{n}(\mathring{\eta})), \mathcal{L}(X_{n}(\eta)) \right) = \mathcal{O}(n^{-d/k(\eta)}) ~.
\end{equation}
Furthermore, the difference between the expected numbers of copies of local configurations $\mathring{\eta}$ and $\eta$ in the graph $G_{n}$ tends to $0$:
\begin{equation}
\label{esp_propre}
\lim_{n\to +\infty} | \EE_{a,b} \lbrack X_{n}(\mathring{\eta}) \rbrack - \EE_{a,b} \lbrack X_{n}(\eta) \rbrack | = 0 ~.
\end{equation}
}
\end{lemm}

In this paper, if $f(n)$ and $g(n)$ are two positive functions, notation $f(n)=\mathcal{O}(g(n))$ means that there exists a constant $C>0$ such that, for all $n$, $f(n)\leq Cg(n)$.\\

\begin{proof}
On the one hand, let us introduce the subset $\mathcal{D}_{r+1}^{\eta}$ of $\mathcal{D}_{r+1}$ defined by:
$$
\mathcal{D}_{r+1}^{\eta} = \left\lbrace \zeta \in \mathcal{D}_{r+1} , \; \forall x \in B(0,r) , \; \zeta(x) = \eta(x)   \right\rbrace ~.
$$
The local configuration $\mathring{\eta}$ is the only element of $\mathcal{D}_{r+1}^{\eta}$ satisfying $k(\mathring{\eta})=k(\eta)$, all the others having at least $k(\eta)+1$ positive vertices. Moreover, the sum of all copies of elements of $\mathcal{D}_{r+1}^{\eta}$ is equal to the number of copies of $\eta$:
\begin{equation}
\label{egal_nb_copies}
X_{n}(\eta) = X_{n}(\mathring{\eta}) + \sum_{\zeta \in \mathcal{D}_{r+1}^{\eta} \setminus \{\mathring{\eta}\}} X_{n}(\zeta) ~.
\end{equation}
On the other hand, the total variation distance between two probability distributions can be written as
\begin{equation}
\label{dTV1}
d_{TV}(\mu, \nu) = \inf \{ \PP(X \not= Y) , \; \mathcal{L}(X) = \mu \; \mbox{ and } \; \mathcal{L}(Y) = \nu \} ~.
\end{equation}
Using this characterization and the identity (\ref{egal_nb_copies}), it follows:
\begin{eqnarray*}
d_{TV} \left( \mathcal{L}(X_{n}(\mathring{\eta})), \mathcal{L}(X_{n}(\eta)) \right) & \leq & \mu_{a,b}( X_{n}(\mathring{\eta}) \not= X_{n}(\eta) ) \\
& \leq & \mu_{a,b}( X_{n}(\eta) > X_{n}(\mathring{\eta}) ) \\
& \leq & \mu_{a,b} \left( \sum_{\zeta \in \mathcal{D}_{r+1}^{\eta} \setminus \{\mathring{\eta}\}} X_{n}(\zeta) > 0 \right) ~.
\end{eqnarray*}
The above sum is an integer valued variable. So, its probability of being positive is bounded by its expectation. Hence,
\begin{eqnarray*}
d_{TV} \left( \mathcal{L}(X_{n}(\mathring{\eta})), \mathcal{L}(X_{n}(\eta)) \right) & \leq & \sum_{\zeta \in \mathcal{D}_{r+1}^{\eta} \setminus \{\mathring{\eta}\}} \EE_{a,b} \lbrack X_{n}(\zeta) \rbrack \\
& \leq & |\mathcal{D}_{r+1}^{\eta}| M_{1} e^{2a(n)} ~,
\end{eqnarray*}
by Lemma \ref{lemm:fonda}. Thus, using $e^{2a(n)}=cn^{-d/k(\eta)}$, we deduce that the total variation distance between the distributions of $X_{n}(\mathring{\eta})$ and $X_{n}(\eta)$ is a $\mathcal{O}(n^{-d/k(\eta)})$.\\
Finally, the same is true for the absolute value of the difference between the expectations of $X_{n}(\mathring{\eta})$ and $X_{n}(\eta)$ since:
$$
| \EE_{a,b} \lbrack X_{n}(\mathring{\eta}) \rbrack - \EE_{a,b} \lbrack X_{n}(\eta) \rbrack | = \sum_{\zeta \in \mathcal{D}_{r+1}^{\eta} \setminus \{\mathring{\eta}\}} \EE_{a,b} \lbrack X_{n}(\zeta) \rbrack ~.
$$
\end{proof}

If the random variable $X_{n}(\mathring{\eta})$ converges weakly as $n$ tends to infinity to a limit $\nu$, inequality (\ref{dvt_propre}) implies that the same is true for $X_{n}(\eta)$. Consequently, replacing $r$ with $r+1$ and without loss of generality, we can assume that vertices of the reference ball $B(0,r)$ which are at distance $r$ from the center $0$, all belong to $V_{-}(\eta)$ (as in Figure \ref{fig:local_configuration}), i.e. $\eta$ can be assumed clean.\\

Now, the geometry (in the sense of the graph structure) of the set $V_{+}(\eta)$ of positive vertices of $\eta$ takes place in our study. Precisely, let us define the \textit{perimeter} $\gamma(\eta)$ of a local configuration $\eta\in\mathcal{D}_{r}$ by the formula:
$$
\gamma(\eta) = \mathcal{V}\; | V_{+}(\eta) | - 2\; | \{ \{x,y\} \in V_{+}(\eta) \times V_{+}(\eta), \; \{x,y\} \in E_{n} \} | ~,
$$
where $\mathcal{V}$ represents the number of neighbors of a vertex. In particular, the perimeter of a local configuration is always an even integer. For instance, that of Figure \ref{fig:local_configuration} is equal to $58$. If $\eta$ is clean, its perimeter $\gamma(\eta)$ merely becomes:
$$
\gamma(\eta) = | \{ \{x,y\} \in V_{+}(\eta) \times V_{-}(\eta)\;, \; \{x,y\} \in E_{n} \} | ~.
$$
In this case, $\gamma(\eta)$ represents the number of pairs of neighboring vertices $x$ and $y$ of $B(0,r)$ having opposite spins under $\eta$. As we shall see in the proof of Lemma \ref{lemm:proba_cond}, the perimeter of a clean local configuration easily occurs in the expression of its local energy. It is the reason why we reduce our study to that of clean local configurations.\\
The following lemma will play an essential role in the proofs of Theorems \ref{theo:PA} and \ref{theo:TV}: it gives a uniform bound for the random variable $n^{d} \mu_{a,b}(I_{x}^{\eta}=1 | \mathcal{F}(\delta\!B(x,r)))$.

\begin{lemm}
\label{lemm:proba_cond}
\textit{
Let $\eta$ be a clean local configuration of radius $r$ and suppose that $e^{2a(n)}=cn^{-d/k(\eta)}$, for some constant $c>0$. Then, there exists a constant $M_{2}>0$ such that for all $n$, for all vertex $x\in V_{n}$ and for all configuration $\sigma\in\mathcal{X}_{n}$,
\begin{equation}
\label{ineg_proba_cond}
c^{k(\eta)} e^{-2 b \gamma(\eta)} (1 - M_{2} e^{2a(n)}) \leq n^{d} \mu_{a,b}(I_{x}^{\eta}=1 | \sigma_{\delta B(x,r)}) \leq c^{k(\eta)} e^{-2 b \gamma(\eta)} ~,
\end{equation}
and for all $n$,
\begin{equation}
\label{ineg_esperance}
c^{k(\eta)} e^{-2 b \gamma(\eta)} (1 - M_{2} e^{2a(n)}) \leq \EE_{a,b} \lbrack X_{n}(\eta) \rbrack \leq c^{k(\eta)} e^{-2 b \gamma(\eta)} ~.
\end{equation}
}
\end{lemm}

Since the quantity $e^{2a(n)}$ tends to $0$ as $n$ tends to infinity, the inequalities (\ref{ineg_proba_cond}) and (\ref{ineg_esperance}) yield the two following limits. For any vertex $x$ and any configuration $\sigma$,
\begin{equation}
\label{lim_proba_cond_esperance}
\lim_{n \to +\infty} n^{d} \mu_{a,b}(I_{x}^{\eta}=1 | \sigma_{\delta B(x,r)}) = \lim_{n \to +\infty} \EE_{a,b} \lbrack X_{n}(\eta) \rbrack = c^{k(\eta)} e^{-2 b \gamma(\eta)} ~.
\end{equation}
Thanks to (\ref{esp_propre}) of Lemma \ref{lemm:propre}, the latter limit is valid for any element of $\mathcal{D}_{r}$ (not necessary clean).

\begin{proof}
Let $x$ be a vertex of $V_{n}$ and denote by $B$ the ball $B(x,r)$. Since the expectations of the variables $X_{n}(\eta)$ and $n^{d} \mu_{a,b}(I_{x}^{\eta}=1 | \mathcal{F}(\delta\!B))$ are equal (see (\ref{egal_expectation})), (\ref{ineg_esperance}) is an immediate consequence of (\ref{ineg_proba_cond}). So, let us prove this relation.\\
Let us start with inserting the perimeter $\gamma(\eta)$ in the expression of the local energy of $\eta$. Assume that $\eta_{x}$ occurs on $B$. Then, there are $\gamma(\eta)$ edges $\{y,\!z\}\in E_{n}$ with $y,z\in B$ satisfying $\eta_{x}(y)\eta_{x}(z)=-1$ and $\alpha(r)-\gamma(\eta)$ ones satisfying $\eta_{x}(y)\eta_{x}(z)=1$. Hence, for all $\sigma\in\mathcal{X}_{n}$, the local energy $H^{B}(\eta_{x}\sigma_{\delta B})$ can be expressed as:
$$
a(n) (2 k(\eta)-\beta(r)) + b \left( \alpha(r)-2\gamma(\eta) + \sum_{{\scriptstyle \{ y,z \} \in E_{n}} \atop \scriptstyle y \in B, z \in \delta B} (-1) \sigma_{\delta B}(z) \right) ~.
$$
The factor $(-1)$ in the latter sum comes from the fact that, by hypothesis, vertices at distance $r$ from $x$ are all negative. Let $\eta'\in\mathcal{D}_{r}$ be a local configuration of radius $r$ with $k(\eta')$ positive vertices. Then, the difference $H^{B}(\eta_{x}'\sigma_{\delta B})-H^{B}(\eta_{x}\sigma_{\delta B})$ between the local energies of $\eta_{x}'$ and $\eta_{x}$ is equal to:
$$
H^{B}(\eta_{x}' \sigma_{\delta B})-H^{B}(\eta_{x} \sigma_{\delta B}) = 2a(n) (k(\eta') - k(\eta)) + b \left( 2\gamma(\eta) + \mathcal{Q}(\eta_{x}') \right) ~,
$$
where
$$
\mathcal{Q}(\eta_{x}') = \sum_{{\scriptstyle \{ y,z \} \in E_{n}} \atop \scriptstyle y,z \in B} \eta_{x}'(y) \eta_{x}'(z) - \alpha(r) + \sum_{{\scriptstyle \{ y,z \} \in E_{n}} \atop \scriptstyle y \in B, z \in \delta B} (\eta_{x}'(y) + 1) \sigma_{\delta B}(z) ~.
$$
The real $\mathcal{Q}(\eta_{x}')$ does not depend on $n$ and can be bounded uniformly on configurations $\eta_{x}'$, $\sigma$: $|\mathcal{Q}(\eta_{x}')| \leq 2\alpha(r+1)$. Moreover, note that the null local configuration of radius $r$ satisfies $\mathcal{Q}(\eta_{x}^{0})=0$. So, using (\ref{ordre_a(n)}), the quantity $H^{B}(\eta_{x}^{0}\sigma_{\delta B})-H^{B}(\eta_{x} \sigma_{\delta B})$ becomes:
\begin{eqnarray*}
H^{B}(\eta_{x}^{0} \sigma_{\delta B}) - H^{B}(\eta_{x} \sigma_{\delta B}) & = & -2a(n) k(\eta) + 2b \gamma(\eta) \\
& = & \log \left( \frac{n^{d}}{c^{k(\eta)}} \right)  + 2b \gamma(\eta) ~,
\end{eqnarray*}
for any configuration $\sigma \in \mathcal{X}_{n}$. Then, using the explicit formula for the conditional probability $\mu_{a,b}(I_{x}^{\eta}=1 | \sigma_{\delta B})$ (relation (\ref{proba_cond})), we get:
\begin{eqnarray*}
\mu_{a,b}(I_{x}^{\eta}=1 | \sigma_{\delta B}) & \leq & e^{H^{B}(\eta_{x} \sigma_{\delta B}) - H^{B}(\eta_{x}^{0} \sigma_{\delta B})} \\
& \leq & \frac{c^{k(\eta)} e^{-2 b \gamma(\eta)}}{n^{d}} ~,
\end{eqnarray*}
i.e. the upper bound of (\ref{ineg_proba_cond}). The lower bound is obtained as follows. For any configuration $\sigma\in\mathcal{X}_{n}$:
\begin{eqnarray*}
n^{d} \mu_{a,b}(I_{x}^{\eta}=1 | \sigma_{\delta B}) & = & \frac{n^{d}}{\sum_{\eta' \in \mathcal{D}_{r}} e^{H^{B}(\eta_{x}' \sigma_{\delta B}) - H^{B}(\eta_{x} \sigma_{\delta B})}} \\
& = & \frac{c^{k(\eta)} e^{-2 b \gamma(\eta)}}{1 + \sum_{\eta', k(\eta')>0} e^{2a(n) k(\eta')+b\mathcal{Q}(\eta_{x}')}} \\
& \geq & \frac{c^{k(\eta)} e^{-2 b \gamma(\eta)}}{1 + e^{2a(n)} \mid \mathcal{D}_{r} \mid e^{2|b|\alpha(r+1)}} ~.
\end{eqnarray*}
Let $M_{2}$ denote the quantity $|\mathcal{D}_{r}| e^{2|b|\alpha(r+1)}$; it only depends on the pair potential $b$ and the radius $r$. Finally, the inequality
$$
\forall u > -1, \; \frac{1}{1 + u} \geq 1 - u
$$
implies the lower bound of (\ref{ineg_proba_cond}).
\end{proof}

\section{Poisson approximation}
\label{section:PA}

This section is devoted to the proof of Theorem \ref{theo:PA}; the distribution of $X_{n}(\eta)$ converges weakly to the Poisson distribution with parameter $c^{k(\eta)} e^{-2 b \gamma(\eta)}$. Previous notations and hypotheses still hold. In particular, the relation (\ref{ordre_a(n)}) between the magnetic field $a(n)$ and the number $k(\eta)$ of positive vertices in $\eta$.\\
Before proving this result, it is worth pointing out here the role of the pair potential $b$. First, remark that local configurations of radius $r$ having the same number of positive vertices can have different perimeters. Theorem \ref{theo:PA} assures that the probability (for $\mu_{a,b}$) of the local configuration $\eta$ of occurring in the graph is asymptotically equal to
$$
1 - e^{-c^{k(\eta)} e^{-2b \gamma(\eta)}} ~.
$$
So, if $b>0$ (resp. $b<0$), this asymptotic probability is a decreasing (resp. increasing) function of the perimeter $\gamma(\eta)$. In other words, if $b>0$ (resp. $b<0$), among the local configurations having the same number of positive vertices, those having the highest asymptotic probability of occurring in the infinite graph are those having the smallest (resp. largest) perimeter.\\
If the pair potential $b$ is null then the perimeter $\gamma(\eta)$ of the local configuration $\eta$ has no influence. All local configurations having the same number of positive vertices have the same asymptotic probability $1-e^{-c^{k(\eta)}}$ of occurring in the graph. In the $2$-dimensional case, this is Theorem $2.4$ of \cite{CoupierDesolneuxYcart2}.\\

In order to prove Theorem \ref{theo:PA}, we use the moment method based on the following lemma (\cite{Bollobas} p.~25).

\begin{lemm}
\label{lemm:moment}
\textit{
Let $(Y_{n})_{n \in \mathbb{N}^{\ast}}$ be a sequence of integer valued, nonnegative random variables and $\lambda$ be a strictly positive real. For all $n , l \in \mathbb{N}^{\ast}$ define $M_{l}(Y_{n})$, the $l$\textnormal{-th moment} of $Y_{n}$, by
\begin{eqnarray*}
M_{l}(Y_{n}) & = & \EE \lbrack Y_{n} ( Y_{n} - 1 ) \ldots ( Y_{n} - l + 1 ) \rbrack \\
& = & \sum_{k \geq l} \PP(Y_{n} = k)\frac{k!}{(k-l)!} ~.
\end{eqnarray*}
If, for all $l \in \mathbb{N}^{\ast}$, $\lim_{n \to \infty} M_{l}(Y_{n}) = \lambda^{l}$ then the distribution of $Y_{n}$ converges weakly as $n$ tends to infinity to the Poisson distribution with parameter $\lambda$.
}
\end{lemm}

First, note that Lemma \ref{lemm:propre} reduces the proof of Theorem \ref{theo:PA} to a clean local configuration $\eta$.\\
So as to lighten formulas, the quantity $c^{k(\eta)} e^{-2b \gamma(\eta)}$ will be simply denoted by $\lambda$. Thanks to Lemma \ref{lemm:moment}, we just need to prove the convergence of $M_{l}(X_{n}(\eta))$ to $\lambda^{l}$ for every $l \in \mathbb{N}^{\ast}$. The case $l=1$ has already been treated at the end of the previous section, where it was proved that $M_{1}(X_{n}(\eta)) = \EE_{a,b} \lbrack X_{n}(\eta) \rbrack$ tends to $\lambda$. From now on, fix an integer $l \geq 2$. In our case, the variable $X_{n}(\eta)$ counts the number of copies in the graph $G_{n}$ of the local configuration $\eta$. Then, the quantity $M_{l}(X_{n}(\eta))$ can be interpreted as the expected number of ordered $l$-tuples of copies of $\eta$.\\
If $\mathcal{B}$ represents a set of balls of radius $r$ whose centers belong to $V_{n}$, then two elements $B(y,r)$ and $B(z,r)$ of $\mathcal{B}$ will be said to be \textit{connected} if there exists an integer $m$ and balls $B_{1},\ldots,B_{m} \in \mathcal{B}$ such that $B_{1}=B(y,r)$, $B_{m}=B(z,r)$ and for $j=1,\ldots,m-1$, $\overline{B}_{j}\,\cap\,B_{j+1} \not= \emptyset$. This last condition allows balls $B_{j}$ and $B_{j+1}$ to be disjoint but their centers are at distance from each other at most $2r+1$. The \textit{connectivity} is an equivalence relation on the set $\mathcal{B}$.\\
For $s=1,\ldots,l$, denote by $\mathcal{C}_{l}(s)$ the set of $l$-tuples of vertices $(x_{1},\ldots,x_{l})$ such that the set of balls $\{ B(x_{1},r),\ldots,B(x_{l},r) \}$ is composed of $s$ equivalence classes for the connectivity relation. Then, the $l$-th moment of $X_{n}(\eta)$ becomes:
$$
M_{l}(X_{n}(\eta)) = \sum_{s=1}^{l} \EE_{a,b} \left\lbrack \sum_{(x_{1},\ldots,x_{l}) \in \mathcal{C}_{l}(s)} I_{x_{1}}^{\eta} \times\ldots\times I_{x_{l}}^{\eta} \right\rbrack ~.
$$
The term corresponding to $s=l$ in the above sum will be denoted by $M_{l}'(X_{n}(\eta))$ and the remaining sum by $M_{l}^{''}(X_{n}(\eta))$. We are going to prove the two following limits
\begin{equation}
\label{Ml'}
\lim_{n \to \infty} M_{l}'(X_{n}(\eta)) = \lambda^{l} ~,
\end{equation}
\begin{equation}
\label{Ml''}
\lim_{n \to \infty} M_{l}^{''}(X_{n}(\eta)) = 0 ~,
\end{equation}
from which Theorem \ref{theo:PA} follows.\\
Let us first check the cardinality of $\mathcal{C}_{l}(s)$.

\begin{lemm}
\label{lemm:cardinal}
\textit{
There exists a constant $C>0$ such that
\begin{equation}
\label{Cls}
\forall s = 1,\ldots ,l , \; | \mathcal{C}_{l}(s) | \leq C n^{ds} ~.
\end{equation}
Furthermore, the cardinality of $\mathcal{C}_{l}(l)$ is equivalent to $n^{dl}$:
\begin{equation}
\label{Cll}
\lim_{n\to +\infty} \frac{| \mathcal{C}_{l}(l) |}{n^{dl}} = 1 ~.
\end{equation}
}
\end{lemm}

\begin{proof}
Let $(x_{1},\ldots,x_{m})$, $m\leq l$, be a $m$-tuple of vertices such that the set of balls $\{B(x_{1},r),\ldots,B(x_{m},r)\}$ is composed of only one equivalence class for the connectivity relation. Each vertex $x_{j}$, $j=1,\ldots,m$, necessary belongs to the ball $B(x_{1},R)$ with $R=l(2r+1)$. So, the number of such $m$-tuples $(x_{1},\ldots,x_{m})$ is bounded above by $\beta(R)^{l-1} n^{d}$. Applying this argument to each of the $s$ equivalence classes of the set formed of the $l$ balls centered at vertices of a given $l$-tuple $(x_{1},\ldots,x_{l})\in\mathcal{C}_{l}(s)$ provides:
$$
| \mathcal{C}_{l}(s) | \leq  \beta(R)^{s(l-1)} n^{ds} ~.
$$
Note that $C=\beta(R)^{l(l-1)}$ only depends on integers $l$ and $r$ and satisfies (\ref{Cls}).\\
Now, we want to choose $l$ vertices $x_{1},\ldots,x_{l}$ such that $(x_{1},\ldots,x_{l})\in\mathcal{C}_{l}(l)$. For the first vertex $x_{1}$, there are $n^{d}$ possibilities. Let $2 \leq j \leq l$ and suppose vertices $x_{1},\ldots,x_{j-1}$ have been chosen. For the $j-$th choice, the set of all vertices $x$ such that $dist(x,x_{k}) \leq 2(r+1)$ for some $1\leq k\leq j-1$, must be avoided. The cardinality of this set is bounded above by $(j-1)\times\,\beta(2(r+1))$ whatever $x_{1},\ldots,x_{j-1}$. This bound does not depend on $n$. As a consequence, from the inequalities
$$
\prod_{1\leq j\leq l} \left( n^{d} - (j-1) \beta(2r+1) \right) \leq | \mathcal{C}_{l}(l) | \leq n^{dl} ~,
$$
estimate (\ref{Cll}) follows.
\end{proof}

Let us prove the first limit (\ref{Ml'}). Let $(x_{1}, \ldots , x_{l}) \in \mathcal{C}_{l}(l)$. By definition of the connectivity relation, note that, for $1\leq i,j\leq l$ and $i\not= j$, no vertex of the ball $B(x_{i},r)$ can be a neighbor of a vertex of the ball $B(x_{j},r)$. The Gibbs measure $\mu_{a,b}$ yields a Markov random field with respect to neighborhoods defined in (\ref{def:voisinage}) (see for example \cite{Malyshev}, Lemma 3 p.~7). As a consequence,
\begin{eqnarray}
\label{prod_proba}
\mu_{a,b} \left( \prod_{i=1}^{l} I_{x_{i}}^{\eta}=1 \right) & = & \EE_{a,b} \left\lbrack \mu_{a,b} \left( \left. \prod_{i=1}^{l} I_{x_{i}}^{\eta}=1 \; \right| \mathcal{F}(\cup_{j=1}^{l} \delta B(x_{j},r)) \right) \right\rbrack \nonumber \\
& = & \EE_{a,b} \left\lbrack \prod_{i=1}^{l} \mu_{a,b} \left( \left. I_{x_{i}}^{\eta}=1 \; \right| \mathcal{F}(\cup_{j=1}^{l} \delta B(x_{j},r)) \right) \right\rbrack \nonumber \\
& = & \EE_{a,b} \left\lbrack \prod_{i=1}^{l} \mu_{a,b} \left( \left. I_{x_{i}}^{\eta}=1 \; \right| \mathcal{F}(\delta B(x_{i},r)) \right) \right\rbrack ~.
\end{eqnarray}
Since the local configuration $\eta$ is clean, Lemma \ref{lemm:proba_cond} and relation (\ref{prod_proba}) provide a control of the probability $\mu_{a,b}(I_{x_{1}}^{\eta}\times\ldots\times I_{x_{l}}^{\eta}=1)$:
$$
\frac{\lambda^{l}}{n^{dl}} \left( 1 - M_{2} e^{2a(n)} \right)^{l} \leq \mu_{a,b} \left( \prod_{i=1}^{l} I_{x_{i}}^{\eta}=1 \right) \leq \frac{\lambda^{l}}{n^{dl}} ~,
$$
uniformly on the $l$-tuple $(x_{1}, \ldots , x_{l}) \in \mathcal{C}_{l}(l)$. Hence,
$$
\frac{| \mathcal{C}_{l}(l) |}{n^{dl}} \lambda^{l} \left( 1 - M_{2} e^{2a(n)} \right)^{l} \leq M_{l}'(X_{n}(\eta)) \leq \frac{| \mathcal{C}_{l}(l) |}{n^{dl}} \lambda^{l} ~.
$$
Finally, as the ratio $|\mathcal{C}_{l}(l)|$ divided by $n^{dl}$ tends to $1$ (relation (\ref{Cll})), the quantity $M_{l}'(X_{n}(\eta))$ tends to the searched limit.\\

There remains to prove that $M_{l}^{''}(X_{n}(\eta))$ tends to $0$ as $n$ tends to infinity. The intuition is that if the local configuration $\eta$ occurs on two balls $B(x,r)$ and $B(x',r)$ with $dist(x,x') \leq 2(r+1)$, then locally (strictly) more than $k(\eta)$ positive vertices are present in a ball of radius $2(r+1)$. This has vanishing probability by Lemma \ref{lemm:fonda}.\\
Let us prove that every term of the sum defining $M_{l}^{''}(X_{n}(\eta))$ tends to $0$: fix an integer $1 \leq s \leq l-1$. Let $(x_{1},\ldots,x_{l})\in\mathcal{C}_{l}(s)$. The set of balls with radius $r$, centered at these vertices, splits into $s$ equivalence classes, say $EC(1),\ldots,EC(s)$. Let us denote by $C_{j}$ the union of balls belonging to the equivalence class $EC(j)$. Once again, we use the markovian character of the Gibbs measure $\mu_{a,b}$:
\begin{eqnarray*}
\mu_{a,b} \left( \prod_{i=1}^{l} I_{x_{i}}^{\eta}=1 \right) & = & \EE_{a,b} \left\lbrack \mu_{a,b} \left( \left. \prod_{i=1}^{l} I_{x_{i}}^{\eta}=1 \; \right| \mathcal{F}(\cup_{j=1}^{l} \delta C_{j}) \right) \right\rbrack \\
& = & \EE_{a,b} \left\lbrack \prod_{j=1}^{s} \mu_{a,b} \left( \left. \prod_{i,B(x_{i},r) \in EC(j)} I_{x_{i}}^{\eta}=1 \; \right| \mathcal{F}(\delta C_{j}) \right) \right\rbrack .
\end{eqnarray*}
As a consequence of $s \leq l-1$, there exists at least one connected component, say $EC(1)$, having at least two elements: let $x(1)$ and $x'(1)$ be two vertices satisfying $dist(x(1),x'(1))\leq 2r+1$ and $B(x(1),r)$, $B(x'(1),r)\in EC(1)$. For every $j=2,\ldots,s$, denote by $x(j)$ one of centers of balls belonging to $EC(j)$. Then, we can write:
$$
\mu_{a,b} \left( \left. \prod_{i,B(x_{i},r) \in EC(j)} I_{x_{i}}^{\eta}=1 \; \right| \mathcal{F}(\delta C_{j})\right) \hspace*{4cm}
$$
\begin{eqnarray*}
\hspace*{3cm} & \leq & \mu_{a,b} \left( \left. I_{x(j)}^{\eta}=1 \; \right| \mathcal{F}(\delta C_{j})\right)\\
& \leq & \EE_{a,b} \left\lbrack \left. \EE_{a,b} \left\lbrack \left. I_{x(j)}^{\eta} \; \right| \mathcal{F}(\delta C_{j} \cup \delta B(x(j),r)) \right\rbrack \; \right| \mathcal{F}(\delta C_{j}) \right\rbrack\\
& \leq & \EE_{a,b} \left\lbrack \left. \EE_{a,b} \left\lbrack \left. I_{x(j)}^{\eta} \; \right| \mathcal{F}(\delta B(x(j),r)) \right\rbrack \; \right| \mathcal{F}(\delta C_{j}) \right\rbrack\\
& \leq & \frac{\lambda}{n^{d}}
\end{eqnarray*}
by Lemma \ref{lemm:proba_cond} ($\eta$ is clean). This last inequality allows us to write:
\begin{eqnarray}
\label{ineg_EC1}
\mu_{a,b} \left( \prod_{i=1}^{l} I_{x_{i}}^{\eta}=1 \right) & \leq & \left( \frac{\lambda}{n^{d}} \right)^{s-1} \EE_{a,b}\left\lbrack \mu_{a,b} \left( \left. \prod_{{\scriptstyle i,B(x_{i},r)} \atop \scriptstyle \in EC(1)} I_{x_{i}}^{\eta}=1 \; \right| \mathcal{F}(\delta C_{1}) \right) \right\rbrack \nonumber \\
& \leq & \left( \frac{\lambda}{n^{d}} \right)^{s-1} \mu_{a,b}\left( \prod_{{\scriptstyle i,B(x_{i},r)} \atop \scriptstyle \in EC(1)} I_{x_{i}}^{\eta}=1 \right) \nonumber \\
& \leq & \left( \frac{\lambda}{n^{d}} \right)^{s-1} \mu_{a,b}( I_{x(1)}^{\eta} = I_{x'(1)}^{\eta} = 1 ) ~.
\end{eqnarray}
Denote by $\mathcal{D}_{2r+1}^{>k(\eta)}$ the (finite) set of local configurations of radius $2r+1$ having at least $k(\eta)+1$ positive vertices. Then, the event $I_{x(1)}^{\eta}=I_{x'(1)}^{\eta}=1$ implies that one of the elements of $\mathcal{D}_{2r+1}^{>k(\eta)}$ occurs in $B(x(1),2r+1)$. It follows that:
\begin{eqnarray*}
\mu_{a,b}( I_{x(1)}^{\eta} = I_{x'(1)}^{\eta} = 1 ) & \leq & \sum_{\zeta \in \mathcal{D}_{2r+1}^{>k(\eta)}} \EE_{a,b} \lbrack I_{x(1)}^{\zeta} \rbrack \\
& \leq & n^{-d} \sum_{\zeta \in \mathcal{D}_{2r+1}^{>k(\eta)}} \EE_{a,b} \lbrack X_{n}(\zeta) \rbrack \\
& \leq & n^{-d} \; | \mathcal{D}_{2r+1}^{>k(\eta)} | \; M_{1} \; e^{2 a(n)} ~,
\end{eqnarray*}
by Lemma \ref{lemm:fonda}. As a consequence, the following bound does not depend on the $l$-tuple $(x_{1},\ldots,x_{l})\in\mathcal{C}_{l}(s)$:
$$
\mu_{a,b} \left( \prod_{i=1}^{l} I_{x_{i}}^{\eta}=1 \right) \leq \frac{\lambda^{s-1}}{n^{ds}} \; | \mathcal{D}_{2r+1}^{>k(\eta)} | \; M_{1} \; e^{2 a(n)} ~.
$$
Finally, Lemma \ref{lemm:cardinal} implies that:
$$
\EE_{a,b} \left\lbrack \sum_{(x_{1}, \ldots , x_{l}) \in \mathcal{C}_{l}(s)} I_{x_{1}}^{\eta} \times \ldots \times I_{x_{l}}^{\eta} \right\rbrack \leq C \lambda^{s-1} \; | \mathcal{D}_{2r+1}^{>k(\eta)} | \; M_{1} \; e^{2 a(n)} ~,
$$
which tends to $0$ as $n$ tends to infinity. Theorem \ref{theo:PA} follows.

\section{The ferromagnetic case}
\label{section:ferromagnetic}

In this section, we suppose the pair potential $b$ is nonnegative and the magnetic field $a(n)$ satisfies the relation (\ref{ordre_a(n)}):
$$
e^{2a(n)} = c n^{-d/k(\eta)} ~,
$$
for some positive constant $c>0$. Under these hypotheses, Theorem \ref{theo:TV} says that the total variation distance between $\mathcal{L}(X_{n}(\eta))$ and its Poisson approximation $\mathcal{P}(c^{k(\eta)} e^{-2b \gamma(\eta)})$ is a $\mathcal{O}(n^{-d/k(\eta)})$.\\
We believe that $n^{-d/k(\eta)}$ is the real speed at which the total variation distance between $\mathcal{L}(X_{n}(\eta))$ and $\mathcal{P}(c^{k(\eta)} e^{-2b \gamma(\eta)})$ tends to zero. Indeed, it seems to be true for the upper bound given by Lemma \ref{lemm:chenstein} (for more details, see Chapter 3 of \cite{Barbour}). Besides, in the case where the local configuration $\eta$ represents a single positive vertex (with $k(\eta)=1$, $\gamma(\eta)=4$, $\rho=1$ and relative to $\| \cdot \|_{1}$), Ganesh et al. \cite{GHOSU} proved that
$$
\frac{ \log d_{TV}(\mathcal{L}(X_{n}(\eta)), \mathcal{P}(c e^{-8b}))}{\log n^{-d}} \; \rightarrow \; 1 ~,
$$
as $n$ tend to infinity.\\

The rest of this section is devoted to the proof of Theorem \ref{theo:TV}. First, note that Lemma \ref{lemm:propre} reduces the proof to a clean local configuration $\eta$.\\
Let us start with some notations and definitions. There is a natural partial ordering on the configuration set $\mathcal{X}_{n}=\{-1,+1\}^{V_{n}}$ defined by $\sigma \leq \sigma'$ if $\sigma(x) \leq \sigma'(x)$ for all vertices $x \in V_{n}$. A function $f:\mathcal{X}_{n} \to \RR$ is \textit{increasing} if $f(\sigma) \leq f(\sigma')$ whenever $\sigma \leq \sigma'$.\\
From the local configuration $\eta$, let us define the subset $\mathcal{D}_{r}(\eta)$ of $\mathcal{D}_{r}$ by:
$$
\mathcal{D}_{r}(\eta) = \{ \eta' \in \mathcal{D}_{r}, V_{+}(\eta') \supset V_{+}(\eta) \} ~.
$$
Each local configuration of $\mathcal{D}_{r}^{\ast}(\eta) = \mathcal{D}_{r}(\eta)\setminus \{\eta\}$ has at least $k(\eta)+1$ positive vertices. Moreover, by definition of $\mathcal{D}_{r}(\eta)$ and for all $x \in V_{n}$, the indicator $\overline{I}_{x}^{\eta}$ defined by
$$
\overline{I}_{x}^{\eta} = \sum_{\eta' \in \mathcal{D}_{r}(\eta)} I_{x}^{\eta'}
$$
is an increasing function. Let us introduce the corresponding random variable $\overline{X}_{n}(\eta)$:
\begin{eqnarray}
\label{sum_increasing}
\overline{X}_{n}(\eta) & = & \sum_{x \in V_{n}} \overline{I}_{x}^{\eta} \\
\label{Xbar}
& = & X_{n}(\eta) + \sum_{\eta' \in \mathcal{D}_{r}^{\ast}(\eta)} X_{n}(\eta') ~,
\end{eqnarray}
whose expectation $\EE_{a,b} \lbrack \overline{X}_{n}(\eta) \rbrack$ will be simply denoted by $\lambda_{n}$. As in the previous Section, the quantity $c^{k(\eta)} e^{-2b \gamma(\eta)}$ will be denoted by $\lambda$.\\
The proof of Theorem \ref{theo:TV} is organized as follows. The total variation distance between $\mathcal{L}(X_{n}(\eta))$ and $\mathcal{P}(\lambda)$ is bounded by:
$$
d_{TV}(\mathcal{L}(X_{n}(\eta)), \mathcal{L}(\overline{X}_{n}(\eta))) + d_{TV}(\mathcal{L}(\overline{X}_{n}(\eta)), \mathcal{P}(\lambda_{n})) + d_{TV}(\mathcal{P}(\lambda_{n}), \mathcal{P}(\lambda)) ~.
$$
Let us respectively denote by $T_{1}$, $T_{2}$ and $T_{3}$ the three terms of the above sum. We are going to prove that each of them is of order $\mathcal{O}(n^{-d/k(\eta)})$. Terms $T_{1}$ and $T_{3}$ are respectively dealt with using Lemmas \ref{lemm:T1} and \ref{lemm:T3}. Applied to the family of indicators $\{\overline{I}_{x}^{\eta}, x \in V_{n}\}$, the Stein-Chen method gives an upper bound for $T_{2}$ (Lemma \ref{lemm:chenstein}). Finally, Lemma \ref{lemm:M2} implies that this upper bound is a $\mathcal{O}(n^{-d/k(\eta)})$.\\


Hypothesis (\ref{ordre_a(n)}) implies that occurrences of local configurations with (strictly) more than $k(\eta)$ positive vertices have vanishing probability. Hence, the random variables $X_{n}(\eta)$ and $\overline{X}_{n}(\eta)$ will be asymptotically equal. So do their expectations.

\begin{lemm}
\label{lemm:T1}
\textit{
The total variation distance between the distributions of $X_{n}(\eta)$ and $\overline{X}_{n}(\eta)$ satisfies:
\begin{equation}
\label{OT1}
d_{TV} \left( \mathcal{L}(X_{n}(\eta)), \mathcal{L}(\overline{X}_{n}(\eta)) \right) = \mathcal{O}(n^{-d/k(\eta)}) ~.
\end{equation}
Furthermore, there exists a constant $M_{3}>0$ such that for all $n$:
\begin{equation}
\label{ineg_lambda_n}
\lambda (1 - M_{3} e^{2 a(n)}) \leq \lambda_{n} \leq \lambda (1 + M_{3} e^{2 a(n)}) ~.
\end{equation}
}
\end{lemm}

\begin{proof}
Thanks to relation between $X_{n}(\eta)$ and $\overline{X}_{n}(\eta)$ (\ref{Xbar}) and characterization (\ref{dTV1}) of the total variation distance, we get:
\begin{eqnarray*}
d_{TV} \left( \mathcal{L}(X_{n}(\eta)), \mathcal{L}(\overline{X}_{n}(\eta)) \right) & \leq & \mu_{a,b}(X_{n}(\eta) \not= \overline{X}_{n}(\eta)) \\
& \leq & \mu_{a,b} \left( \sum_{\eta' \in \mathcal{D}_{r}^{\ast}(\eta)} X_{n}(\eta') > 0 \right) \\
& \leq & \sum_{\eta' \in \mathcal{D}_{r}^{\ast}(\eta)} \EE_{a,b} \lbrack X_{n}(\eta') \rbrack \\
& \leq & | \mathcal{D}_{r}^{\ast}(\eta) | \; M_{1} e^{2 a(n)} ~,
\end{eqnarray*}
by Lemma \ref{lemm:fonda}. So, the quantity $d_{TV}(\mathcal{L}(X_{n}(\eta)), \mathcal{L}(\overline{X}_{n}(\eta)))$ is equal to $\mathcal{O}(e^{2 a(n)})= \mathcal{O}(n^{-d/k(\eta)})$.\\
Using the previous inequalities, a control of the expectation $\lambda_{n}$ of $\overline{X}_{n}(\eta)$ is obtained:
$$
\EE_{a,b} \lbrack X_{n}(\eta) \rbrack \leq \lambda_{n} \leq \EE_{a,b} \lbrack X_{n}(\eta) \rbrack + | \mathcal{D}_{r}^{\ast}(\eta) | M_{1} e^{2 a(n)} ~.
$$
Since $\eta$ is clean, relation (\ref{ineg_esperance}) can be applied. The above control becomes:
$$
\lambda \left( 1 - M_{2} e^{2 a(n)} \right) \leq \lambda_{n} \leq \lambda \left( 1 + | \mathcal{D}_{r}^{\ast}(\eta) | M_{1} \lambda^{-1} e^{2 a(n)} \right) ~.
$$
Recall that constants $M_{1}$ and $M_{2}$ do not depend on the size $n$. Relation (\ref{ineg_lambda_n}) follows by letting:
$$
M_{3} = \max \{ M_{2} , | \mathcal{D}_{r}^{\ast}(\eta) | M_{1} \lambda^{-1} \} ~.
$$
\end{proof}

Using (\ref{ineg_lambda_n}), we shall now bound $T_{3}$.

\begin{lemm}
\label{lemm:T3}
\textit{
The total variation distance between the Poisson distributions with parameters $\lambda_{n}$ and $\lambda$ satisfies:
$$
T_{3} = d_{TV}(\mathcal{P}(\lambda_{n}), \mathcal{P}(\lambda)) = \mathcal{O}(n^{-d/k(\eta)}) ~.
$$
}
\end{lemm}

\begin{proof}
The total variation distance between two probability distributions on the set of integers can be expressed as:
$$
d_{TV}(\mu, \nu) = \frac{1}{2} \sum_{m \geq 1} | \mu(m) - \nu(m) | ~.
$$
Let $m \geq 1$. Thanks to relation (\ref{ineg_lambda_n}), the difference $\lambda_{n}^{m} e^{-\lambda_{n}} - \lambda^{m} e^{-\lambda}$ is easily controlled. Thus, $d_{TV}(\mathcal{P}(\lambda_{n}), \mathcal{P}(\lambda))$ is bounded by $\frac{1}{2} \sum_{m\geq1} \max \{ \alpha_{m}, \beta_{m} \}$ where:
$$
\alpha_{m} = \frac{\lambda^{m} e^{-\lambda}}{m!} \left( e^{\lambda M_{3} e^{2a(n)}} (1 + M_{3} e^{2a(n)})^{m} - 1 \right)
$$
and
$$
\beta_{m} = \frac{\lambda^{m} e^{-\lambda}}{m!} \left( 1 - e^{-\lambda M_{3} e^{2a(n)}} (1 - M_{3} e^{2a(n)})^{m} \right) ~.
$$
Using the convexity of the function
$$
f_{m}:\; \rbrack -1, 1 \lbrack \longrightarrow \RR , \; x \longmapsto (1 + x)^{m} e^{\lambda x} ~,
$$
one easily checks that $\alpha_{m} \geq \beta_{m}$, for all $m \in \NN^{\ast}$. As a consequence,
$$
d_{TV}(\mathcal{P}(\lambda_{n}), \mathcal{P}(\lambda)) \leq \frac{1}{2} \sum_{m\geq1} \alpha_{m} \leq \frac{1}{2} \left( e^{2 \lambda M_{3} e^{2a(n)}} - 1 \right)
$$
which is of order $\mathcal{O}(n^{-d/k(\eta)})$ by relation (\ref{ordre_a(n)}).
\end{proof}

There remains to bound the term $T_{2} = d_{TV}(\mathcal{L}(\overline{X}_{n}(\eta)), \mathcal{P}(\lambda_{n}))$. This is based on the Stein-Chen method and particulary on Corollary 2.C.4, p.~26 of \cite{Barbour} which is described below (Proposition \ref{prop:chenstein}). Let $\{I_{i}\}_{i\in I}$ be a family of random indicators with expectations $\pi_{i}$. Let us denote
$$
W = \sum_{i\in I} I_{i} \; \mbox{ and } \; \theta = \sum_{i\in I} \pi_{i} ~.
$$
The random variables $\{I_{i}\}_{i\in I}$ are \textit{positively related} if for each $i$, there exists random variables $\{J_{j,i}\}_{j\in I}$ defined on the same probability space such that
$$
\mathcal{L}(J_{j,i} , \; j \in I) = \mathcal{L}(I_{j} , \; j \in I | I_{i} = 1) 
$$
and, for all $j\not= i$, $J_{j,i}\geq I_{j}$.

\begin{prop}
\label{prop:chenstein}
\textit{
If the random variables $\{I_{i}\}_{i\in I}$ are positively related then:
$$
d_{TV}(\mathcal{L}(W), \mathcal{P}(\theta)) \leq \frac{1 - e^{-\theta}}{\theta} \left( Var(W) - \theta + 2 \sum_{i\in I} \pi_{i}^{2} \right) ~.
$$
}
\end{prop}

Proposition \ref{prop:chenstein} can be applied to our context. Indeed, for a positive value of the pair potential $b$, the Gibbs measure $\mu_{a,b}$ defined by (\ref{def:mes_gibbs}) satisfies the FKG inequality, i.e.
\begin{equation}
\label{ineg_FKG}
\EE_{a,b} \lbrack fg \rbrack \geq \EE_{a,b} \lbrack f \rbrack \EE_{a,b} \lbrack g \rbrack ~,
\end{equation}
for all increasing functions $f$ and $g$ on $\mathcal{X}_{n}$: see for instance Section 3 of \cite{FKG}. Then, Theorem 2.G p.~29 of \cite{Barbour} implies that the increasing random indicators $\{\overline{I}_{x}^{\eta}, x\in V_{n}\}$ are positively related. Replacing $I_{i}$ with $\overline{I}_{x}^{\eta}$, $W$ with $\overline{X}_{n}(\eta)$ and $\theta$ with $\lambda_{n}$, Proposition \ref{prop:chenstein} produces the following result. This is the only place where the hypothesis $b\geq 0$ is actually used in the proof.

\begin{lemm}
\label{lemm:chenstein}
\textit{
If the pair potential $b$ is nonnegative then the following inequality holds:
\begin{equation}
\label{boundT2}
d_{TV}(\mathcal{L}(\overline{X}_{n}(\eta)), \mathcal{P}(\lambda_{n})) \leq \frac{1}{\lambda_{n}} \left( Var_{a,b} \lbrack \overline{X}_{n}(\eta)\rbrack - \lambda_{n} + 2 \sum_{x \in V_{n}} \EE_{a,b} \lbrack \overline{I}_{x}^{\eta} \rbrack ^{2} \right) ~.
\end{equation}
}
\end{lemm}

The bound (\ref{boundT2}) indicates that, as $n \to +\infty$, the sum $\sum_{x \in V_{n}} \EE_{a,b} \lbrack \overline{I}_{x}^{\eta} \rbrack ^{2}$ is small and the distance to the Poisson approximation is essentially the difference between the variance and the expectation of $\overline{X}_{n}(\eta)$. Using good estimates on the first two moments of the random variable $\overline{X}_{n}(\eta)$, this difference will be bounded. The case of the first moment of $\overline{X}_{n}(\eta)$, i.e. its expectation $\lambda_{n}$, has been treated in Lemma \ref{lemm:T1}. The following result concerns its second moment:

\begin{lemm}
\label{lemm:M2}
\textit{
The second moment $M_{2}(\overline{X}_{n}(\eta)) = \EE_{a,b}\lbrack \overline{X}_{n}(\eta) (\overline{X}_{n}(\eta) - 1) \rbrack$ of the random variable $\overline{X}_{n}(\eta)$ satisfies:
$$
M_{2}(\overline{X}_{n}(\eta)) = \lambda^{2} + \mathcal{O}(n^{-d/k(\eta)}) ~.
$$
}
\end{lemm}

Writing the variance of the variable $\overline{X}_{n}(\eta)$ as $M_{2}(\overline{X}_{n}(\eta))+\lambda_{n}-\lambda_{n}^{2}$ and the sum $\sum_{x \in V_{n}} \EE_{a,b} \lbrack \overline{I}_{x}^{\eta} \rbrack ^{2}$ as the ratio $\lambda_{n}^{2}/n^{d}$, we deduce from Lemma \ref{lemm:chenstein} that:
$$
T_{2} \leq \frac{1}{\lambda_{n}} \left( M_{2}(\overline{X}_{n}(\eta)) + \lambda_{n}^{2} ( \frac{2}{n^{d}} - 1 ) \right) ~.
$$
The inequalities given by (\ref{ineg_lambda_n}) and Lemma \ref{lemm:M2} allow us to control the expectation $\lambda_{n}$ and the second moment $M_{2}(\overline{X}_{n}(\eta))$ of the random variable $\overline{X}_{n}(\eta)$. This implies:
$$
T_{2} = \mathcal{O}(e^{2 a(n)}) + \frac{2 \lambda^{2}}{n^{d}} ~,
$$
which is a $\mathcal{O}(e^{2 a(n)}) = \mathcal{O}(n^{-d/k(\eta)})$ since $k(\eta) \geq 1$.\\
Let us finish the proof of Theorem \ref{theo:TV} by proving Lemma \ref{lemm:M2}.

\begin{proofbis}{Lemma \ref{lemm:M2}}
First, recall that $\mathcal{C}_{2}(s)$, for $s=1,2$, represents the set of couples $(x_{1},x_{2})$ whose set $\{ B(x_{1},r),B(x_{2},r)\}$ splits into $s$ equivalence classes for the connectivity relation. In other words, $(x_{1},x_{2})$ belongs to $\mathcal{C}_{2}(1)$ if $dist(x_{1},x_{2})\leq 2r+1$ and to $\mathcal{C}_{2}(2)$ otherwise. So, the second moment of $\overline{X}_{n}(\eta)$ is equal to:
$$
M_{2}(\overline{X}_{n}(\eta)) = \sum_{s=1}^{2} \EE_{a,b} \left\lbrack \sum_{(x_{1},x_{2}) \in \mathcal{C}_{2}(s)} \overline{I}_{x_{1}}^{\eta} \times \overline{I}_{x_{2}}^{\eta} \right\rbrack ~.
$$
Each indicator $\overline{I}_{x}^{\eta}$ is defined as the sum of $I_{x}^{\eta'}$, $\eta'\in\mathcal{D}_{r}(\eta)$. Hence, the second moment $M_{2}(\overline{X}_{n}(\eta))$ becomes:
\begin{equation}
\label{dec_Xbar}
M_{2}(\overline{X}_{n}(\eta)) = \sum_{\eta_{1}, \eta_{2} \in \mathcal{D}_{r}(\eta)} \left( E^{1}(\eta_{1}, \eta_{2}) + E^{2}(\eta_{1}, \eta_{2}) \right) ~,
\end{equation}
where for $s=1,2$, the quantity $E^{s}(\eta_{1},\eta_{2})$ is defined by:
$$
E^{s}(\eta_{1}, \eta_{2}) = \EE_{a,b} \left\lbrack \sum_{(x_{1},x_{2}) \in \mathcal{C}_{2}(s)} I_{x_{1}}^{\eta_{1}} I_{x_{2}}^{\eta_{2}} \right\rbrack ~.
$$
Let $(\eta_{1},\eta_{2})$ be a couple of local configurations belonging to $\mathcal{D}_{r}(\eta)$ and $(x_{1},x_{2})$ be a couple of vertices. In a first time, consider $(x_{1},x_{2})\in\mathcal{C}_{2}(1)$. Then it has been already seen at the end of the previous section that the event $I_{x_{1}}^{\eta_{1}}=I_{x_{2}}^{\eta_{2}}=1$ implies that one of the elements of $\mathcal{D}_{2r+1}^{>k(\eta)}$ necessary occurs in $B(x_{1},2r+1)$. It follows that:
\begin{eqnarray*}
\mu_{a,b}( I_{x_{1}}^{\eta_{1}} = I_{x_{2}}^{\eta_{2}} = 1 ) & \leq & \sum_{\zeta \in \mathcal{D}_{2r+1}^{>k(\eta)}} \EE_{a,b} \lbrack I_{x_{1}}^{\zeta} \rbrack \\
& \leq & n^{-d} \sum_{\zeta \in \mathcal{D}_{2r+1}^{>k(\eta)}} \EE_{a,b} \lbrack X_{n}(\zeta) \rbrack \\
& \leq & n^{-d} \; | \mathcal{D}_{2r+1}^{>k(\eta)} | \; M_{1} \; e^{2 a(n)} ~,
\end{eqnarray*}
by Lemma \ref{lemm:fonda}. Thus, thanks to Lemma \ref{lemm:cardinal}, we deduce that $E^{s}(\eta_{1},\eta_{2})$, for $s=1,2$, is a $\mathcal{O}(e^{2 a(n)})=\mathcal{O}(n^{-d/k(\eta)})$.\\
Now, let us suppose that $(x_{1},x_{2})\in\mathcal{C}_{2}(2)$. Some technics already used in the previous section give:
\begin{equation}
\label{prob_egal_prod}
\EE_{a,b} \left\lbrack I_{x_{1}}^{\eta_{1}} I_{x_{2}}^{\eta_{2}} \right\rbrack = \EE_{a,b} \left\lbrack \prod_{i=1}^{2} \mu_{a,b} \left( \left. I_{x_{i}}^{\eta_{i}} = 1 \; \right| \mathcal{F}(\delta B(x_{i},r)) \right) \right\rbrack ~.
\end{equation}
At this point of the proof, two cases must be distinguished: either both local configurations $\eta_{1}$ and $\eta_{2}$ are equal to $\eta$ or not. In the first case, Lemma \ref{lemm:proba_cond} and (\ref{prob_egal_prod}) imply:
$$
\EE_{a,b} \left\lbrack I_{x_{1}}^{\eta} I_{x_{2}}^{\eta} \right\rbrack \leq n^{-2d} \lambda ^{2} ~.
$$
Then, the quantity $E^{2}(\eta,\eta)$ which is actually the second moment of $X_{n}(\eta)$, is bounded by $\lambda^{2}$. In the other case, at least one of the two local configurations $\eta_{1}$, $\eta_{2}\in\mathcal{D}_{r}(\eta)$ is different from $\eta$, i.e. has at least $k(\eta)+1$ positive vertices. Then, coupling Lemma \ref{lemm:fonda} with (\ref{prob_egal_prod}), it follows that:
$$
E^{2}(\eta_{1},\eta_{2}) \leq M_{1}^{2} e^{2 a(n)(k(\eta_{1}) + k(\eta_{2}) - 2 k(\eta))} ~,
$$
which is a $\mathcal{O}(e^{2 a(n)}) = \mathcal{O}(n^{-d/k(\eta)})$.
\end{proofbis}

\bibliographystyle{plain}
\bibliography{pa_ising_bibli}

\end{document}